\documentclass[11pt]{article}
\usepackage{amsfonts}
\usepackage{amsmath}
\usepackage{url}

\newcommand{\B}{\mathbb}

\newtheorem{conjecture}{Conjecture}[section]

\newtheorem{corollary}[conjecture]{Corollary}
\newtheorem{definition}[conjecture]{Definition}

\newtheorem{lemma}[conjecture]{Lemma}

\newtheorem{proposition}[conjecture]{Proposition}

\newtheorem{theorem}[conjecture]{Theorem}

\newcommand{\R}{\mathbb{R}}

\newcommand{\N}{\mathbb{N}}

\newcommand{\teichmuller}{Teichm\"uller }

\newcommand{\mcg}{${\rm{MCG}}(\Sigma)$ }
\newcommand{\pml}{${\rm{PML}}(\Sigma)$ }

\newenvironment{proof}{\noindent {\bf Proof} \hspace{.1cm}}{ $\square$  \vspace{.15cm}}

\title{Free subgroups of surface mapping class groups}
\author{James W. Anderson, Javier Aramayona, Kenneth J. Shackleton\footnote{Partially supported by a Japan Society 
for the Promotion of Science post-doctoral fellowship for foreign researchers, number PE$05043$}}
\date{15 May 2006}

\begin{document}

\maketitle

\begin{abstract}
\noindent We quantify the generation of free subgroups of surface mapping class groups by pseudo-Anosov mapping classes in terms
of their translation distance and the distance between their axes. Our methods make reference to \teichmuller space only.

\medskip

\noindent MSC 20F65 (primary),  57M50 (secondary)
\end{abstract}

\section{Introduction}

Free subgroups of mapping class groups have attracted considerable
attention from a number of authors over a number of years. It is a classical
result, proven independently by Ivanov \cite{ivanov} and by McCarthy \cite{mccarthy}, 
that any collection $A$ of pseudo-Anosov mapping classes with pairwise distinct axes freely generates a free group of rank $|A|$, 
so long as each element is first raised to a sufficiently high power.  The supporting
argument is based on the so-called Ping-Pong Lemma, quoted here in Section 4, and makes reference to Thurston's boundary of \teichmuller space.

The purpose of this work is to establish a ping-pong argument entirely inside \teichmuller space.
This allows us to not only recover the result of Ivanov and of McCarthy by somewhat different means, but also to use the Teichm\"uller metric 
to quantify sufficiently high powers in a natural way. Our work relies on an application of Minsky's Bounded Projection Theorem \cite{minsky}. One may also use Thurston's train tracks to quantify sufficiently high powers, see Hamidi-Tehrani \cite{h-t}.

The plan of this paper is as follows. In Section 2 we introduce
all the background and notation we shall need.
In Section 3 we recover the known result that axes of
independent pseudo-Anosov mapping classes cannot be asymptotic. In Section 4 we exhibit ping-pong sets inside \teichmuller space for independent
pseudo-Anosovs and use these to show that sufficiently high powers of
independent pseudo-Anosov mapping classes generate a free group. In Section 5 we give a quantitative version of this result, 
with Theorem \ref{quantitative-schottky} our aim. Section 4 and Section 5 are logically independent, but both rest firmly on Section 2 and Section 3.

\noindent{\bf{Acknowledgments.}} The second author would like to thank Vlad Markovic, 
Brian Bowditch and Caroline Series for many interesting and helpful conversations. The third author 
wishes to thank both the Japan Society for the Promotion of Science, for its partial support, and Osaka University, 
for its warm hospitality. He also wishes to thank the Institut des Hautes \'Etudes Scientifiques for providing a stimulating working environment.

\section{Background}

We refer the reader to \cite{bonahon}, \cite{imayoshi} and \cite{ivanov-notes} for detailed studies of geodesic and measured geodesic laminations, \teichmuller spaces and mapping class groups, respectively, and recall only what we need. Throughout this paper, a {\em surface} $\Sigma$ will mean an orientable connected surface of negative Euler characteristic, with genus $g$ and $p$ punctures, and with empty boundary.  A {\em curve} on $\Sigma$ is the free homotopy class of a simple closed loop that is neither homotopic to a point nor to a puncture, and we denote by ${\cal S}(\Sigma)$ the set of all curves on $\Sigma$.  We say that two curves are {\em disjoint} if they can be realised disjointly, and if they are not disjoint we say they {\em intersect essentially} or are {\em transverse}.

\subsection{\teichmuller space}

The {\em \teichmuller space} $T(\Sigma)$ is the space of all marked finite
area hyperbolic structures on $\Sigma$, up to homotopy. More
specifically, a point in $T(\Sigma)$ is an equivalence class
$[(\sigma,f)]$, where $\sigma$ is a finite area hyperbolic structure on
$\Sigma$ and $f:\Sigma \longrightarrow \sigma$ is a homeomorphism, called
a {\em marking of} $\sigma$. One declares two pairs $(\sigma,f)$ and
$(\sigma',g)$ to be equivalent if and only if $g \circ f^{-1}$ is
homotopic to an isometry between $\sigma$ and $\sigma'$. To simplify
our notation, we will use the same symbol to denote a
point in $T(\Sigma)$ as to denote a particular marked finite area
hyperbolic structure on $\Sigma$.

The space $T(\Sigma)$ is homeomorphic to $\mathbb{R}^{6g-6+2p}$, and can be compactified by attaching the space \pml of all {\em{projective measured laminations}} on $\Sigma$, topologically a sphere of dimension $6g-7+2p$. The closed ball $T(\Sigma) \cup {\rm PML}(\Sigma)$ is sometimes known as the {\em Thurston compactification of \teichmuller space}, and we shall denote it by $\overline{T(\Sigma)}$.

Any curve $\alpha\in {\cal S}(\Sigma)$ induces a length function $\ell_\alpha(x)$ on $T(\Sigma)$, where $\ell_\alpha(x)$ denotes the length of the unique geodesic representative of the class $\alpha$ in the hyperbolic structure $x$ on $\Sigma$. Given $\epsilon >0$, the {\em{$\epsilon$-thick part}} of
\teichmuller space is defined by 
\[ T_{\geq \epsilon}(\Sigma)  = \{x \in T(\Sigma) : \ell_{\alpha}(x) \geq \epsilon \mbox{ for all } \alpha
\in {\cal S}(\Sigma) \}, \]
and the {\em{$\epsilon$-thin part}} of \teichmuller space is defined by
\[ T_{\leq \epsilon}(\Sigma)  = \{x \in T(\Sigma) : \ell_{\alpha}(x) \leq \epsilon \mbox{ for all } \alpha
\in {\cal S}(\Sigma) \}. \]
The space $T(\Sigma)$ admits two natural metrics, the
\teichmuller metric and the Weil-Petersson metric. In this paper
we will only consider the \teichmuller metric, and refer the reader
to \cite{wolpert-wp} for a thorough study of the Weil-Petersson
geometry. Given points $x, y \in T(\Sigma)$, the {\em{\teichmuller
distance}} between $x=[(\sigma,f)]$ and $y=[(\sigma',g)]$ is defined as $${\rm d}_T(x,y)
=  \frac{1}{2} \inf\{ \log(K(h))\},$$ where the infimum ranges
over all the quasiconformal homeomorphisms $h$ in the homotopy
class of $f \circ g^{-1}$ and $K(h)$ is the dilatation of $h$. By a celebrated
result of Teichm\"uller, for any two points in Teichm\"uller space there is a unique quasiconformal homeomorphism $h$ (in the appropriate homotopy class) realising the distance between the points. Endowed with the \teichmuller metric, $T(\Sigma)$ is a uniquely
geodesic, proper metric space. It is worth noting, however, the
\teichmuller metric is not non-positively curved in any standard
sense (see \cite{masur-class}, \cite{masur-wolf},
\cite{mccarthy-actions}, \cite{bowditch-ECM}).

For subsets $X$ and $Y$ of $T(\Sigma)$, we define the {\em nearest point distance} ${\rm d}_T(X,Y)$ to be
\[ {\rm d}_T(X,Y) =\inf\{ {\rm d}_T(x,y)\: :\: x\in X, y\in Y\}. \]

\subsection{Mapping class groups}

The {\em mapping class group} ${\rm MCG}(\Sigma)$ of $\Sigma$ is the group of all homotopy classes of orientation preserving self-homeomorphisms of $\Sigma$. It is a
finitely presented group, generated by a finite collection of Dehn
twists about simple closed curves of $\Sigma$. There is a natural action of \mcg on $T(\Sigma)$ by changing the marking; the
quotient $T(\Sigma) / {\rm MCG}(\Sigma)$ is the {\em moduli space of
$\Sigma$}. Except for a few low-dimensional cases, the
mapping class group \mcg is isomorphic to an index two subgroup of the full isometry group of both the \teichmuller metric and the Weil-Petersson metric, by results of Royden
\cite{royden} and Masur-Wolf \cite{masur-wolf} respectively.

There are some similarities between the action
of \mcg on $T(\Sigma)$, equipped with the \teichmuller metric, and
that of a geometrically finite group acting isometrically on a simply-connected Hadamard manifold of pinched negative curvature. For example, the elements of \mcg can be
classified in a manner that mimics the classification of the
isometries of the pinched Hadamard manifold according to their
dynamics on its ideal boundary. An infinite
order mapping class is either
{\em{reducible}}, so it fixes some non-empty and finite collection of disjoint curves, or is otherwise pseudo-Anosov.

A {\em pseudo-Anosov} mapping class $\phi$ is represented by a pseudo-Anosov diffeomorphism $\Sigma$. That is, there exists a real number $r = r(\phi) >1$, the dilatation of $\phi$, such that for any hyperbolic metric $(\Sigma, \sigma)$ there exists a unique diffeomorphism $f$ representing $\phi$ and two measured laminations, $\lambda^-$ and $\lambda^+$, geodesic in $\sigma$ such that $f(\lambda^+) = r\lambda^+$ and
$f(\lambda^-) = \frac{1}{r} \lambda^-$. We shall write $\lambda^{\pm}$ to denote either element of the set $\{\lambda^+, \lambda^-\}$. The measured lamination $\lambda^{\pm}$ satisfies the following three fundamental
properties \cite{thurston}:

\begin{enumerate} 

\item $\lambda^{\pm}$ is {\em{uniquely ergodic}}: if $\mu$ is a measured lamination whose support ${\rm{supp}}(\mu)$ is equal to ${\rm{supp}}(\lambda^{\pm})$, then
$\mu$ and $\lambda^{\pm}$ are proportional;

\item $\lambda^{\pm}$ is {\em{minimal}}: if $\mu$ is a measured lamination satisfying ${\rm{supp}}(\mu) \subseteq {\rm{supp}}(\lambda^{\pm})$, 
then either ${\rm{supp}}(\mu) = \emptyset$ or ${\rm{supp}}(\mu) = {\rm{supp}}(\lambda^{\pm})$, and

\item $\lambda^{\pm}$ is {\em{maximal}}: if $\mu$ is a measured lamination satisfying ${\rm{supp}}(\lambda^{\pm}) \subseteq {\rm{supp}}(\mu)$, then ${\rm{supp}}(\mu) = {\rm{supp}}(\lambda^{\pm})$.

\end{enumerate}

We will say a projective measured lamination is {\em uniquely ergodic} if one (and hence any) representative of its projective class is a uniquely ergodic measured lamination, and whenever $\lambda$ is uniquely ergodic we shall also use, where there can be no ambiguity, $\lambda$ to denote both a measured lamination and its projective class. 

The fixed point set ${\rm{Fix}}(\phi)$ of $\phi$ in $\overline{T(\Sigma)}$ is precisely $\{\lambda^+, \lambda^-\}$. These fixed points behave like attracting and repelling fixed points 
for $\phi$. More specifically, with $s \in \{ +1, -1\}$,
for any neighbourhood $U$ of $\lambda^s$ in $\overline{T(\Sigma)}$ and any compact set $K$ in $\overline{T(\Sigma)}\setminus \{\lambda^{-s}\}$ we have $\phi^{sn}(K) \subseteq U$ for sufficiently large $n$ 
(see \cite{ivanov-notes}). It is known that a pseudo-Anosov mapping
class $\phi$ fixes a bi-infinite \teichmuller geodesic, the
{\em{axis} of} $\phi$, on which it acts by translation. By the above discussion, the set of accumulation points
of this axis on \pml is ${\rm{Fix}}(\phi) = \{\lambda^+, \lambda^-\}$.

The {\em{translation distance}} ${\rm Tr}(\phi) = {\rm{inf}} \{{\rm d}_T(x,\phi(x))
: x \in T(\Sigma) \}$ of a pseudo-Anosov mapping class is
always realised, and is always realised on the axis of
$\phi$. Furthermore, both ${\rm Tr}$ and the property of being pseudo-Anosov are invariant under conjugation. The following result is due to Ivanov \cite{ivanov-anosov}.

\begin{theorem}[Ivanov \cite{ivanov-anosov}]
For a surface $\Sigma$ and $L>0$, there are only finitely many conjugacy classes of pseudo-Anosov mapping classes of translation distance at most $L$.\label{finite-anosov}
\end{theorem}

It follows there exists a constant $\ell_{\rm{min}} = \ell_{\rm{min}}(\Sigma)  > 0$ such that all pseudo-Anosov mapping classes in \mcg have translation distance at least $\ell_{\rm{min}}$. Lower bounds for $\ell_{\rm{min}}$, in terms of the topological type of $\Sigma$, 
have been found by Penner \cite{penner}.

The following terminology is due to Minsky \cite{minsky}.

\begin{definition} For a surface $\Sigma$ and $\epsilon >0$, a
\teichmuller geodesic $c$ is said to be {\em $\epsilon$-precompact} if $c$ is entirely
contained in the $\epsilon$-thick part $T_{\geq \epsilon}(\Sigma)$ of
$T(\Sigma)$.
\end{definition}

For any pseudo-Anosov mapping class $\phi$, the projection of its axis into moduli space is compact. Moreover, by the continuity of the length functions $x \longrightarrow \ell_\alpha(x)$, there is a uniform lower bound on all $\ell_\alpha(x)$, where $\alpha \in
{\cal S}(\Sigma)$ and $x$ lies on the axis of $\phi$. Applying Theorem \ref{finite-anosov} to remove all dependence on $\phi$ yields the following.

\begin{corollary}
For a surface $\Sigma$ and $L>0$, there exists a positive real number $\epsilon = \epsilon(L, \Sigma)$ such that, if $\phi$
is any pseudo-Anosov mapping class of translation distance at most $L$, the axis of $\phi$ is $\epsilon$-precompact.
\label{thick-axes}
\end{corollary}

Let $\phi, \psi \in {\rm MCG}(\Sigma)$ be two pseudo-Anosov mapping
classes and let ${\rm{Fix}}(\phi)$ and ${\rm{Fix}}(\psi)$ be their respective fixed point sets in $\overline{T(\Sigma)}$. It is known (see
\cite{mccarthy-papadopoulos}) that either ${\rm{Fix}}(\phi)=
{\rm{Fix}}(\psi)$, in which case $\phi$ and $\psi$ are powers of the same
pseudo-Anosov mapping class, or ${\rm{Fix}}(\phi) \cap {\rm{Fix}}(\psi) =
\emptyset$. This prompted the following definition.

\begin{definition} [\cite{mccarthy-papadopoulos}]
For a surface $\Sigma$, two pseudo-Anosov mapping classes $\phi, \psi \in
{\rm MCG}(\Sigma)$ are said to be {\rm independent} if ${\rm{Fix}}(\phi) \cap
{\rm{Fix}}(\psi) = \emptyset$.
\end{definition}

\section{Divergence of pseudo-Anosov axes}

It is one consequence of Minsky's Bounded Projection Theorem
\cite{minsky}, found by Farb-Mosher \cite{farb-mosher}, that the axes of two
independent pseudo-Anosov mapping classes cannot be asymptotic in
$T(\Sigma)$, that is their Hausdorff distance is not finite. 
We offer a proof of a very much related result, namely that the distance function between the axes of independent pseudo-Anosov mapping classes is a proper function. (The corresponding result for 
the Weil-Petersson metric is due Daskalopoulos-Wentworth \cite{daskalopoulos-wentworth}.) We shall not need the
 Bounded Projection Theorem at all, but instead a theorem of Wolpert \cite{wolpert} and a special case of Lemma 2.1 from \cite{dkw}.

\begin{theorem} [\cite{wolpert}]
For a surface $\Sigma$ and $D>0$, let $x,y \in T(\Sigma)$ with ${\rm d}_T(x,y) \leq D$. Then,
for any curve $\alpha \in {\cal S}(\Sigma)$, we have
$$e^{-2D}\ell_\alpha(x) \leq l_\alpha(y) \leq e^{2D} \ell_\alpha(x).$$
\label{wolpert}
\end{theorem}

Although the intersection number of two arbitrary projective measured laminations is not well-defined, we can still decide whether their intersection number should be zero or non-zero.

\begin{lemma} [\cite{dkw}] For a surface $\Sigma$, let $\lambda, \lambda' \in {\rm{PML}}(\Sigma)$.  Let $(x_n)$ be a sequence of points in
$T(\Sigma)$ converging to $\lambda$ and let
$(\alpha_n)$ be a sequence of curves converging to $\lambda'$ as measured laminations.
Suppose there exists $R>0$ such that $\ell_{\alpha_n}(x_n) \leq R$. Then, $i(\lambda, \lambda') = 0$. 
\label{dkw}
\end{lemma}

We are now ready to prove the distance function restricted to the
axes of a pair of pseudo-Anosov mapping classes is a proper map. Recall, a map is said to be {\em proper} if the preimage of any compact subset of the range is compact.

\begin{proposition} For a surface $\Sigma$, let $\phi, \psi \in {\rm{MCG}}(\Sigma)$ be two
independent pseudo-Anosov mapping classes and let $c,c':\R
\longrightarrow T(\Sigma)$ be arc-length parametrizations of their
respective axes. Then, the map $(t,s) \longrightarrow {\rm d}_T(c(t),c'(s))$ is
a proper map.
\label{proper-map}
\end{proposition}

\begin{proof} Suppose the result is not true. Then, there are unbounded sequences $(t_n)$ and $(s_n)$ in $\R$ and a real number $M>0$ such that ${\rm d}_T(c(t_n),c'(s_n)) \leq M$ for
all $n \in \N$. By passing to subsequences if need be, we may assume that $(t_n)$ and $(s_n)$ are each either monotonically
increasing or monotonically decreasing.  Let us suppose they are both monotonically increasing, as the remaining cases can be treated analogously.

Let $\lambda, \mu \in {\rm{PML}}(\Sigma)$ be such that
$x_n = c(t_n) \longrightarrow \lambda$ and $y_n = c'(s_n) \longrightarrow
\mu$ in $\overline{T(\Sigma)}$, as $n \longrightarrow \infty$.  Then, for all $n\in \N$, there
exists $x'_n \in c(\R)$ such that ${\rm d}_T(x_n,x'_n) \leq {\rm Tr}(\phi)$ and $x'_n
= \phi^{k_n}(x_0)$ for some $x_0 \in c$ and $k_n = k(n)$. In
particular, we have $k_n \longrightarrow \infty$ as $n\longrightarrow \infty$. Since ${\rm d}_T(x_n,y_n) \leq M$, we see that
${\rm d}_T(x'_n,y_n) \leq M+{\rm Tr}(\phi)$. Denote by $M'$ the upper bound $M + {\rm Tr}(\phi)$.

Choose any curve $\alpha \in {\cal S}(\Sigma)$ and let $\alpha_n = \phi^{k_n}(\alpha)$, noting $(\alpha_n)$ converges to $\lambda$ as measured laminations.
It follows from the definition of the action of \mcg on $T(\Sigma)$ that
$\ell_{\alpha_n}(x'_n) = \ell_{\alpha}(x_0)$. This fact, together with
the upper bound in Theorem \ref{wolpert}, implies $\ell_{\alpha_n}(y_n) \leq e^{2M'} \ell_{\alpha}(x_0)$.
 Therefore, the sequences $(y_n)$ and $(\alpha_n)$ satisfy the
hypotheses of Lemma \ref{dkw}, with $R$ equal to $e^{2M'}\ell_{\alpha}(x_0)$, and we conclude
$i(\lambda,\mu)=0$. As $\lambda$ and $\mu$ are both maximal and both minimal, we find ${\rm{supp}}(\lambda) = {\rm{supp}}(\mu)$ and thus
$\lambda = \mu$ since $\lambda$ and $\mu$ are uniquely ergodic. However, this is contrary to the assumption that $\phi$ and $\psi$ be independent as pseudo-Anosov mapping classes.
\end{proof}

An immediate consequence of Proposition \ref{proper-map} is the following.

\begin{corollary}  [\cite{farb-mosher}] For a surface $\Sigma$, let $\phi$ and $\psi$ be two independent pseudo-Anosov mapping classes and let
 $c,c' : \B{R} \longrightarrow T(\Sigma)$ be arc-length parametrizations of their respective axes. Then, the Hausdorff distance between 
 $c(\B{R})$ and $c'(\B{R})$ is infinite.\label{time-interval-length}
\end{corollary}

\section{Ping-ponging in \teichmuller space}

The purpose of this section is to exhibit ping-pong sets inside
\teichmuller space for any given finite family of independent pseudo-Anosov mapping classes. This recovers the theorem of Ivanov \cite{ivanov} and McCarthy \cite{mccarthy}. In Section 5 we shall give a quantitative version of this result.

For a closed subset $C$ of $T(\Sigma)$ and a point $x \in T(\Sigma)$, one defines the closest-point projection of $x$ into $C$ as 
\[ \pi_C(x) = \{ y\in C\: :\: {\rm d}_T(x, y) \leq  {\rm d}_T(x, z) {\mbox{ for all }} z \in C \}. \] 
Note $\pi_C(x)$
is non-empty since $T(\Sigma)$ is a proper space. For another subset $C' \subset T(\Sigma)$, we define
\[ \pi_C(C') =\bigcup_{x \in C'}\pi_C(x). \]
Given any path $c : \B{R} \longrightarrow T(\Sigma)$ we shall also use $\pi_{c}$ to denote $\pi_{c(\B{R})}$, the projection to the image of $c$. The next result is Minsky's
Bounded Projection Theorem, a union of Contraction Theorem (1), Corollary 4.1 and Theorem 4.2 from \cite{minsky}, and highlights some of the hyperbolic
behaviour of any thick part of \teichmuller space. For its statement, we shall first need the following standard definition.

\begin{definition} Let $(X, d_X)$ and $(Y, d_Y)$ be metric spaces, and let $K \geq 1, \kappa \geq 0$. A $(K, \kappa)${\rm -quasi-isometric 
embedding of $X$ into $Y$} is a map $f: X \longrightarrow Y$ such that 

$$\frac{1}{K} d_X(x,x') - \kappa \leq d_Y(f(x),f(x')) \leq K d_X(x,x') + \kappa,$$  

for all $x,x' \in X$. A $(K, \kappa)${\rm-quasi-geodesic in $Y$} is a $(K, \kappa)$-quasi-isometric embedding of a closed subinterval of $\mathbb{R}$ into $Y$. 
\end{definition}

In what follows, $B_{r}(x)$ denotes the ball in $T(\Sigma)$ of radius $r$ and centre $x$.

\begin{theorem} [\cite{minsky}]
For a surface $\Sigma$ and $\epsilon > 0$, there exists a constant $b = b(\epsilon, \Sigma)$ such that the following hold.\begin{enumerate}
\item Given any $\epsilon$-precompact
geodesic $c$ and any $x \in T(\Sigma)$, we have
$${\rm {diam}}(\pi_c(B_{{\rm d}_T(x,c(\B{R}))}(x))) \leq b.$$
\item Given any $\epsilon$-precompact
geodesic $c$ and any $x,y \in T(\Sigma)$, we have $${\rm{diam}}(\pi_c(x) \cup \pi_c(y)) \leq {\rm d}_T(x,y) +4b.$$
\item Given $K\geq 1$ and $\kappa \geq 0$, there exists a non-negative real number $M=M(K,\kappa,\epsilon, \Sigma)$ such that the following holds: If $q$ is
a $(K,\kappa)$-quasi-geodesic path in $T(\Sigma)$ whose endpoints
are connected by an $\epsilon$-precompact \teichmuller geodesic
$c$, then the image of $q$ is contained in the $M$-neighbourhood of the image of $c$.
\end{enumerate}\label{projection}
\end{theorem}

The next result will be a key ingredient for our ping-pong
argument. Roughly speaking, it says precompact \teichmuller
geodesics diverge ``sufficiently fast".
We shall make use of the following notation: Given an embedding 
$c:\R \longrightarrow T(\Sigma)$ and two of its points $x=c(t), x'=c(t')$, we will write $x < x'$ if $t < t'$.

\begin{proposition} For a surface $\Sigma$ and $\epsilon>0$, let $c,c': \R \longrightarrow T(\Sigma)$ be arc-length parametrizaions of two $\epsilon$-precompact \teichmuller geodesics so that $O = c(0)$ and $O' = c'(0)$ realise the nearest point distance $D$ between $c(\B{R})$ and $c'(\B{R})$. Then, there exist two points $P^+, P^- \in c(\B{R})$ and two points 
$Q^+, Q^- \in c'(\B{R})$, with $P^-<O<P^+$ and $Q^-<O'<Q^+$, such that the following hold:

\begin{enumerate}

   \item For all $x \in c(\R)$ with $x>P^+$ and all $y \in c'(\R)$ with
   $y>Q^+$, $${\rm d}_T(x,y) > {\rm{max}} \{{\rm d}_T(O,x), {\rm d}_T(O', y)\}, \mbox{ {\rm{and}}}$$

   \item For all $x \in c(\R)$ with $x<P^-$ and all $y \in c'(\R)$ with
   $y<Q^-$, $${\rm d}_T(x,y) > {\rm{max}} \{{\rm d}_T(O,x), {\rm d}_T(O', y)\}.$$
\end{enumerate}
\label{fast-divergence}
\end{proposition}

\begin{proof}
We show only the first part of the proposition, as the second part follows by an analogous argument. Suppose, for contradiction, the statement
is not true. Then, we can find two unbounded sequences
$(x_n)$ and $(y_n)$ of points in $c(\B{R})$ and $c'(\B{R})$, respectively, such
that $O < x_n < x_{n+1}$, $O'< y_n < y_{n+1}$ and ${\rm d}_T(x_n,y_n) \leq
{\rm{max}}\{{\rm d}_T(O,x_n), {\rm d}_T(O', y_n)\}$ for all $n \in \N$. Passing to a further subsequence if need be, we have ${\rm d}_T(O,x_n)
\geq {\rm d}_T(O',y_n)$ for all $n$, or we have ${\rm d}_T(O,x_n)
\leq {\rm d}_T(O',y_n)$ for all $n$. Without loss of generality, let us suppose the former holds.

Let $g_n$ be the $c(\R)$-segment between $O$ and $x_n$, noting the length of $g_n$
tends to infinity as $n \longrightarrow \infty$. Let $q_n$ be the concatenation of
the unique \teichmuller geodesic segment from $O$ to $O'$, the
$c'(\R)$-segment from $O'$ to $y_n$, and the unique \teichmuller
geodesic segment from $y_n$ to $x_n$. Then,
$${\rm{length}}(q_n) \leq D + {\rm d}_T(O',y_n) + {\rm d}_T(x_n,y_n) \leq D + {\rm d}_T(O, x_{n}) + {\rm d}_T(O, x_{n}) =  D + 2 {\rm{length}}(g_n).$$ Moreover, $${\rm length}(g_{n}) = {\rm d}_{T}(O, x_{n}) \leq {\rm d}_{T}(O, O') + {\rm d}_{T}(O', y_{n}) + {\rm d}_{T}(y_{n}, x_{n}) = {\rm length}(q_{n}).$$
We deduce $${\rm length}(g_{n}) \leq {\rm length}(q_{n}) \leq D + 2{\rm length}(g_{n}),$$ and therefore $q_n$ is a $(2,D)$-quasi-geodesic. In particular,
notice that the constants of quasi-geodesicity are independent of
$n$. By Theorem \ref{projection}(3), there is a constant $M=M(2, D, \epsilon, \Sigma)$ such that the image of $q_n$ is contained in the $M$-neighbourhood of the image of $g_n$, for all $n
\in \N$.  Since ${\rm d}_{T}(O', y_{n})$ tends to infinity as $n \longrightarrow \infty$, this implies that the Hausdorff distance between $c([0, \infty))$ and $c'([0, \infty))$ is at most $M$. According to Proposition \ref{proper-map}, this is a contradiction. \end{proof}

\begin{corollary}  For a surface $\Sigma$ and $\epsilon>0$, let $c,c': \R \longrightarrow T(\Sigma)$ be arc-length parametrizations of two 
$\epsilon$-precompact \teichmuller geodesics. Then, the set
$\pi_{c'}(c(\R))$ is bounded.
\label{bounded-projection}
\end{corollary}

\begin{proof} Let $D$ be the nearest-point distance between $c(\B{R})$ and $c'(\B{R})$, and choose $O \in c(\B{R})$ and $O' \in c'(\B{R})$ such that ${\rm d}_T(O,O') = D$. Note $O' \in \pi_{c'}(O)$ and $O \in
\pi_c(O')$. Suppose, for contradiction, that $\pi_{c'}(c(\R))$ is not bounded.
Then, there exists a sequence $(x_n)$ of points from $c(\B{R})$ with
${\rm d}_T(O', y_n) \longrightarrow \infty$, for any $y_n \in \pi_{c'}(x_n)$.
Furthermore, the sequence $(x_n)$ satisfies ${\rm d}_T(O, x_n) \longrightarrow \infty$, by Theorem \ref{projection}(2). Passing to further subsequences if need be, we may assume $x_n>P^+$ and $y_n>Q^+$, say, where
$P^+ \in c(\B{R})$ and $Q^+ \in c'(\B{R})$ are the respective points given by Proposition \ref{fast-divergence}.

We have ${\rm d}_T(x_n,y_n) > {\rm{max}}\{{\rm d}_T(O,x_n), {\rm d}_T(O',y_n) \}$ and, according to Theorem \ref{projection}(1), we also have ${\rm{diam}}(\pi_{c'}(B_{{\rm d}_T(x_n,c'(\B{R})}(x_n))) \leq b$, for all $n \in \N$. However, $O \in B_{{\rm d}_T(x_n,c'(\B{R}))}(x_n)$ and so $O' \in \pi_{c'}(B_{{\rm d}_T(x_n,c'(\B{R})}(x_n))$ for all $n$. From this we deduce ${\rm d}_T(O', y_{n}) \leq b$ for all $n$, and this is a contradiction. \end{proof}

Given an arc-length parametrization  $c : \B{R} \longrightarrow T(\Sigma)$ of a Teichm\"uller geodesic and a real number
$R>0$, we introduce the subsets 
\[ \Pi(c,R) = \{ x \in T(\Sigma) : \pi_c(x) \subset c([R, \infty)) \} \]
and 
\[ \Pi(c,-R) = \{ x \in
T(\Sigma) : \pi_c(x) \subset c((-\infty, \-R]) \}\]
of $T(\Sigma)$. We note that, if $R < R'$, then $\Pi(c, R') \subset \Pi(c, R)$ and $\Pi(c, -R') \subset \Pi(c, -R)$.

\begin{corollary} For a surface $\Sigma$ and $\epsilon>0$, let $c,c': \R \longrightarrow T(\Sigma)$ be arc-length parametrizations of two $\epsilon$-precompact \teichmuller geodesics. Then, there exists $R>0$ such that the sets $\Pi(c,R)$, $\Pi(c,-R)$, $\Pi(c',R)$ and $\Pi(c',-R)$ are pairwise disjoint. \label{ping-pong-sets}
\end{corollary}

\begin{proof} That $\Pi(c,R) \cap \Pi(c,-R)$
and $\Pi(c',R) \cap \Pi(c',-R)$ are both the empty set for sufficiently large $R$ is a
trivial consequence of Minsky's Bounded Projection Theorem. Let us just show $\Pi(c,R) \cap \Pi(c',R) = \emptyset$,
since the remaining cases can be proven analogously.
We again argue by contradiction, by supposing that for all $n \in \N$
there exists $x_n \in \Pi(c,n) \cap \Pi(c',n)$. Let $y_n \in
\pi_c(x_n)$ and let $z_n \in \pi_{c'}(x_n)$, for all $n$. Note that, in particular, $(y_n)_{n\in\N}$ and $(z_n)_{n\in\N}$ are unbounded
sequences on the geodesics $c(\B{R})$ and $c'(\B{R})$, respectively. On passing to subsequences if need be, we have 
${\rm d}_T(x_n,y_n) \geq {\rm d}_T(x_n,z_n)$, for all $n$, or ${\rm d}_T(x_n,y_n) \leq {\rm d}_T(x_n,z_n)$ for all $n$. 
Without loss of generality, we assume the former holds.

The diameter of the projection of $B_{{\rm d}_T(x_n,y_n)}(x_n)$ into $c(\B{R})$ is at most $b$, by Theorem \ref{projection}(1). In particular, ${\rm d}_T(\pi_c(z_n),\pi_c(x_n))
\leq b$ and it follows ${\rm d}_T(\pi_{c}(c'), y_{n}) \leq b$ for all $n \in \N$. This is a contradiction, since $\pi_{c}(c')$ has bounded diameter, 
by Corollary \ref{bounded-projection}, and $(y_n)$ is an unbounded sequence on the geodesic $c(\B{R})$. \end{proof}

Corollary \ref{ping-pong-sets} gives us enough information to apply the following lemma, the statement of which is recorded from \cite{bridson-haefliger}, and deduce sufficiently high powers of $n$ independent pseudo-Anosov mapping classes freely generate a free group of rank $n$.

\begin{lemma}[Ping-Pong Lemma] Let $X$ be a set and let
$f_1, \dots, f_n$ be bijections from $X$ to itself. Suppose, for every
$i=1,\dots,n$, there exist pairwise disjoint subsets $A_1^+, A_1^-,
\dots, A_n^+, A_n^-$ of $X$ such that $f_i(X\setminus A_i^-)
\subseteq A_i^+$ and $f_i^{-1}(X\setminus A_i^+) \subseteq A_i^-$, for each $i$.
Then, under composition, $f_1, \dots, f_n$ freely generate a free group of rank $n$.\label{pp}
\end{lemma}

\begin{corollary} [\cite{ivanov}, \cite{mccarthy}] For a surface $\Sigma$, let $\phi_1, \dots, \phi_n$ be pairwise independent 
pseudo-Anosov mapping classes in ${\rm MCG}(\Sigma)$. Then, there exists a natural number $N$ such that $\phi_1^N, \dots, \phi^N_n$ freely generate a free group of rank $n$. \label{schottky}
\end{corollary}

\begin{proof}
Let $c_i:\R \longrightarrow T(\Sigma)$ be an arc-length parametrization of the axis of $\phi_i$ for each $i=1, \dots,n$. By Corollary \ref{thick-axes} there exists a real number $\epsilon >0$ such that $c_i$ is $\epsilon$-precompact for each
$i=1, \dots,n$. We note $\epsilon$ depends only on the maximal translation distance among the $\phi_{i}$. 

Now Corollary \ref{ping-pong-sets} applied to all pairs $c_j$ and $c_k$ of parametrizations for $1\le j\ne k\le n$ implies there exists $R>0$
such that the sets $\Pi(c_1,\pm R), \dots,
\Pi(c_n,\pm R)$ are all pairwise disjoint. Let $\ell_{\rm{min}}>0$ be the
the minimal translation distance among all pseudo-Anosov mapping classes in MCG($\Sigma$). Let $N$ be the least integer such that $N>{2R} / {\ell_{\rm{min}}}$. The mapping classes $\phi_1^N, \dots, \phi_n^N$ and the sets
$A_1^{\pm} = \Pi(c_1,\pm R), \dots, A_n^{\pm} = \Pi(c_n,\pm R)$ satisfy the hypotheses of Lemma \ref{pp}, and we conclude $\phi_1^N, \dots, \phi_n^N$ freely generate a free
group of rank $n$.
\end{proof}

\section{A quantitative ping-pong argument for pseudo-Anosovs}

The purpose of this section is to give the promised quantitative version of Corollary \ref{schottky}. On the way, we will show quantitative versions of Proposition \ref{fast-divergence} and Corollary \ref{bounded-projection} for the axes of pseudo-Anosov mapping classes. Let us begin with the following definition.

\begin{definition} For a surface $\Sigma$, $R>0$, and $x \in T(\Sigma)$, we say a curve $\alpha \in {\cal S}(\Sigma)$ is {\rm $R$-short on $x$}
 if $\ell_\alpha(x) \leq R$, and we let $S_R(x)$ be the set of $R$-short curves on $x$.
\end{definition}

The next result is a special case of the main result in \cite{birman-series}, where Birman-Series show the number of simple closed geodesics on a given surface grows at most polynomially in the length bound. 
The degree of this polynomial depends only on the topological type of the surface. We remark an improved version of this result has been given by Rivin \cite{rivin}.

\begin{theorem}[\cite{birman-series}, \cite{rivin}] For a surface $\Sigma$ and $R>0$, 
there exists an integer $B=B(R, \Sigma)$ such that the cardinality of the set $S_R(x)$ is at most $B$ for all $x \in T(\Sigma)$.
\label{bound-on-curves}
\end{theorem} 

Recall that a {\em pants decomposition for $\Sigma$} is a maximal collection of pairwise distinct and disjoint curves on $\Sigma$. It is a theorem of Bers that there exists a universal constant $R_* = R_*(\Sigma)$ such that every point $x \in T(\Sigma)$ has a pants decomposition whose curves each have length at most $R_*$ in $x$.

Let $\epsilon>0$ and consider the thick part $T_{\geq \epsilon}(\Sigma)$ of the
\teichmuller space of $\Sigma$. A simple area argument shows that, given $\epsilon >0$,
there is a constant $R = R(\epsilon)$ such that, for every point $x \in T_{\geq \epsilon}(\Sigma)$, the set $S_R(x)$ contains a 
pants decomposition of $\Sigma$ and a curve transverse to each curve in the pants decomposition. From this discussion, and from Theorem 
\ref{finite-anosov} and Corollary \ref{thick-axes}, we obtain
the following lemma.

\begin{lemma} For a surface $\Sigma$ and $L>0$, there exists a real number $F=F(L, \Sigma)>0$ such that, for  any pseudo-Anosov of translation distance at most $L$ 
and any point $x$ on its axis, the set $S_F(x)$ contains a pants decomposition of $\Sigma$ and a curve intersecting each curve in the pants decomposition essentially.
\label{F-constant}
\end{lemma}

Suppose that $\phi, \psi $ are independent pseudo-Anosov mapping classes of translation distance at most $L>0$, with $c,c':\R \longrightarrow T(\Sigma)$ arc-length parametrizations of their respective axes. By Corollary \ref{thick-axes}, there exists $\epsilon = \epsilon(L, \Sigma)>0$ such that $c(\B{R})$ and $c'(\B{R})$ are $\epsilon$-precompact
\teichmuller geodesics. By Theorem \ref{projection}(3), for any $D>0$ there exists a non-negative real number $M=M(2, D, \epsilon, \Sigma)$ such that any $(2,D)$-quasi-geodesic connecting the ends of an $\epsilon$-precompact geodesic lies in the $M$-neighbourhood of the geodesic. Let $F=F(L, \Sigma)$
be the constant given by Lemma \ref{F-constant}, and recall $B=B(e^{2(M+L)}F, \Sigma) \geq 0$ is a uniform upper bound on the number of $e^{2(M+L)}F$-short curves over all points in $T(\Sigma)$. Note, $B$ depends only on $D, L$ and $\Sigma$. We have the following result, a quantitative analogue of Proposition \ref{fast-divergence}.

\begin{proposition} For a surface $\Sigma$, let $\phi, \psi $ be independent pseudo-Anosov mapping classes of translation distance at most $L>0$. Let $c,c':\R \longrightarrow T(\Sigma)$ 
be arc-length parametrizations of their respective axes so that $O = c(0)$ and $O' = c'(0)$ together realise the nearest-point distance 
$D$ of $c(\B{R})$ and $c'(\B{R})$. For $R= {\rm{max}}\{B!+2, (B!+2)L\}$, the following holds: For all $x \in c(\B{R}) \setminus B_R(O)$ and 
all $y \in c'(\B{R}) \setminus B_R(O')$, 
\begin{center}${\rm d}_T(x,y) >{\rm{max}} \{ {\rm d}_T(O,x), {\rm d}_T(O',y) \}.$
\end{center}\label{quantitative-divergence}
\end{proposition}

\begin{proof}
Suppose the result were not true. Then, there are points $x \in c(\B{R})$
and $y \in c'(\B{R})$ with ${\rm d}_T(x,O)>R$, ${\rm d}_T(y,O')>R$ and ${\rm d}_T(x,y) \leq
{\rm{max}} \{ {\rm d}_T(O,x),{\rm d}_T(O',y) \}$. Assume without loss of generality that ${\rm d}_T(O,x) \geq {\rm d}_T(O',y)$.
Let $q$ be the concatenation of the unique geodesic segment
between $O$ and $O'$, the $c'(\R)$-segment between $O'$ and $y$ and
the unique geodesic segment between $y$ and $x$. Since ${\rm d}_T(x,y)
\leq {\rm{max}} \{{\rm d}_T(O,x),{\rm d}_T(O',y)\}$ we have, as in the proof of
Proposition \ref{fast-divergence}, that $q$ is a
$(2,D)$-quasi-geodesic connecting $O$ to $x$. Therefore the image of $q$ is
entirely contained in the $M$-neighbourhood of the unique geodesic connecting $O$ to $x$, where $M = M(2, D,\epsilon, \Sigma)$ as per Theorem 4.2(3).

Let $I=B!+2$. Since the translation distance of $\psi$ is at most $L$ and $R \geq (B!+2)L=IL$, the points 
$y_1=O', y_2 =\psi(O') \dots, y_I=\psi^{I-1}(O')$ all lie on the geodesic from $O'$ to $y$, which is a geodesic subpath 
of the quasi-geodesic $q$. Since the image of $q$ is contained in the $M$-neighbourhood of the geodesic connecting $O$ to $x$, there are
points $x_1, \ldots, x_I$ in $c(\mathbb{R})$ such that $d(x_i,y_i) \leq M$ for all $i=1,\dots, I$. Since the translation distance of $\phi$ is also at most
$L$, it follows that there are points $z_1, \ldots, z_I$ (not necessarily distinct) such that $d_T(x_i,z_i) \leq L$ and $z_i = \phi^{j(i)}(O)$ for $i=1,
\dots, I$. Therefore $d_T(y_i,z_i) \leq M+L$ for $i=1, \dots, I$. 

By Lemma \ref{F-constant} there exists $F=F(L, \Sigma)>0$ such that the set $S_F(O')$ of short curves in $O'$ contains a set $S$ consisting of a pants decomposition
and a single curve intersecting each curve in the pants decomposition essentially. Then, every element of the set $S$ is $e^{2(M+L)}F$-short in $O$ by Theorem \ref{wolpert}.
Now, if a curve $\alpha$ is $F$-short in $O'$ then $\psi^{i-1}(\alpha)$ is $F$-short in $y_i$ for $i=1, \ldots, I$. Therefore $\phi^{-j(i)} \psi^{i-1}(\alpha)$
is $e^{2(M+L)}F$-short in $O$. In particular, all elements of the set $\phi^{-j(i)} \psi^{i-1}(S)$ are $e^{2(M+L)}F$-short in $O$. Moreover, note that
the set $\phi^{-j(i)} \psi^{i-1}(S)$  also consists of a pants decomposition and a curve intersecting each curve in the pants decomposition transversally, since 
$\phi^{-j(i)} \psi^{i-1}$ is a mapping class. 

Since the set of $e^{2(M+L)}F$-short curves in $O$ has cardinality at most $B$ and $I>B!+1$, there must be some $k \in \{0, \ldots, I-1\}$ such that
$\phi^{-j(k)} \psi^{k-1}(\beta)= \beta$ for all $\beta \in S$. It follows $\phi^{j(k)}$ and $\psi^{k-1}$ share the same action on the set $S(\Sigma)$ of all curves on $\Sigma$. 
Therefore $\phi^{j(k)}$ and $\psi^{k-1}$ are either equal or, for only a few exceptional surfaces, perhaps differ by a 
hyperelliptic involution. Regardless, both share common fixed points in $\overline{T(\Sigma)}$ and this is contrary to their independence.  
\end{proof}

\begin{corollary}  For a surface $\Sigma$ and $L>0$, let $\phi, \psi $ be independent pseudo-Anosov mapping classes of translation distance at most $L>0$. 
Let $c,c':\R \longrightarrow T(\Sigma)$ be arc-length parametrizations of their respective axes so that $O = c(0)$ and 
$O' = c'(0)$ together realise the nearest-point distance of $c(\B{R})$ and $c'(\B{R})$. Let $R$ and $b$ be as per Proposition \ref{quantitative-divergence} and Theorem
\ref{projection}, respectively. Then,
\begin{center}$\pi_c(c') \subseteq B(O, R+4b) \cap c(\B{R}) \mbox{ and }  \pi_{c'}(c) \subseteq B(O', R+4b) \cap c'(\B{R}).$\label{quantitative-projection}
\end{center}
\end{corollary}

\begin{proof}
We need only prove one of inclusions. Suppose, for contradiction, that $\pi_{c'}(c)$ is not entirely
contained in $B(O', R+4b) \cap c'(\B{R})$. Then, there are points $x \in
c(\B{R})$ and $y\in \pi_{c'}(x) \subset c'(\B{R})$ with ${\rm d}_T(O',y)
> R+4b$. By Theorem \ref{projection}(2), we have $${\rm d}_T(O,x) + 4b \geq {\rm{diam}}(\pi_{c}(O) \cup \pi_{c}(x)) \geq {\rm d}_{T}(O', y) > R + 4b$$ and hence 
${\rm d}_T(O, x) > R$. According to Proposition
\ref{quantitative-divergence}, we also have ${\rm d}_T(x,y)>{\rm d}_T(O,x)$. In
particular, $O \in B_{{\rm d}_T(x,y)}(x)$. By Theorem
\ref{projection}(1), ${\rm{diam}}(\pi_{c'}(B_{{\rm d}_T(x,y)}(x))) \leq b$ and so ${\rm d}_{T}(O', y) \leq b$. This is a contradiction, 
and we deduce the corollary. \end{proof}

\begin{corollary} For a surface $\Sigma$ and $L>0$, let $\phi, \psi $ be independent pseudo-Anosov mapping classes of translation distance at most $L>0$. Let
 $c,c':\R \longrightarrow T(\Sigma)$ be arc-length parametrizations of their respective axes so that $O = c(0)$ and 
 $O' = c'(0)$ together realise the nearest-point distance of $c(\B{R})$ and $c'(\B{R})$. Let $R$ and $b$ be as per Proposition \ref{quantitative-divergence} and Theorem
\ref{projection}, respectively. Then, the sets $\Pi(c,R+6b)$, 
 $\Pi(c,-R-6b)$, $\Pi(c',R+6b)$ and $\Pi(c',-R-6b)$ are pairwise
disjoint. \label{quantitative-pingpong}
\end{corollary}

\begin{proof} We prove $\Pi(c,R+6b) \cap \Pi(c',R+6b) =\emptyset$. Again,
suppose that there exists $x\in T(\Sigma)$ such that $x \in
\Pi(c,R+6b) \cap \Pi(c',R+6b)$. Let $y \in \pi_c(x), w \in \pi_{c'}(y)$ and $z \in
\pi_{c'}(x)$. Without loss of generality, we assume that ${\rm d}_T(x,y) \geq {\rm d}_T(x,z)$. By Theorem \ref{projection}(1) we have
${\rm{diam}}(\pi_{c'}(B_{{\rm d}_T(x,y)}(x))) \leq b$, and in particular ${\rm d}_{T}(w, z) \leq b$. On the other hand, appealing to 
Corollary 5.5, we have $${\rm d}_{T}(w, z) \geq {\rm d}_{T}(O', z) - {\rm d}_{T}(O', w) \geq R + 6b - R - 4b = 2b > b$$ and this is a contradiction.
\end{proof}

We are ready to give the promised quantitative version of Corollary \ref{schottky}. Recall, $\ell_{\rm{min}}$ is defined as the 
minimal translation distance among all pseudo-Anosov mapping classes from MCG($\Sigma$).

\begin{theorem} For a surface $\Sigma$ and $L>0$, let $\phi_1, \dots, \phi_n$ be pseudo-Anosov mapping classes of
translation distance at most $L$. Let $R$ and $b$ be as per Proposition \ref{quantitative-divergence} and Theorem
\ref{projection}, respectively. If $N>(2R+12b) / {\ell_{\rm{min}}}$,
then $\phi_1^N, \dots, \phi_n^N$ freely generate a free group of
rank $n$. \label{quantitative-schottky}
\end{theorem}

\begin{proof} Corollary \ref{quantitative-pingpong} implies the sets
$\Pi(c_1,\pm(R+6b)),\dots, \Pi(c_n,\pm(R+6b))$ are ping-pong sets
for $\phi_1^N, \dots, \phi_n^N$. The result now follows from the
Ping-Pong Lemma.
\end{proof}

{\footnotesize James W. Anderson\\School of Mathematics\\University of Southampton\\Southampton SO17 1BJ\\England\\j.w.anderson@soton.ac.uk}

{\footnotesize Javier Aramayona\\Mathematics Institute\\University of Warwick\\Coventry CV4 7AL\\England\\jaram@maths.warwick.ac.uk}

{\footnotesize Kenneth J. Shackleton (corresponding author)\\Institute des Hautes \'Etudes Scientifiques\\Le Bois-Marie\\35 Route de Chartres\\F-91440 Bures sur Yvette\\France\\k.j.shackleton@maths.soton.ac.uk}

\end{document}